\documentclass[12pt, dvips, twoside]{amsart}
\usepackage{amssymb,latexsym,bbm,graphicx,epsfig,epic,eepic,oldgerm}
\usepackage{a4wide}

\theoremstyle{plain}
\newtheorem{theorem}{Theorem}[section]
\newtheorem{lemma}[theorem]{Lemma}
\newtheorem{proposition}[theorem]{Proposition}

\newtheorem*{remark*}{Remark}
\newtheorem*{remarks*}{Remarks}
\newtheorem{remark}[theorem]{Remark}

\newtheorem*{example*}{Example}
\newtheorem*{examples*}{Examples}

\newcommand{\proofend}{\hspace*{\fill} $\Box$\\}
\newcommand{\diam}{\hspace*{\fill} $\Diamond$\\}

\def\s{\smallskip}
\def\m{\medskip}
\def\eps{\varepsilon}

\def\im{\operatorname {im}}

\def\Diffc0{\operatorname{Diff^c_0}}
\def\Symp{\operatorname{Symp}}

\def\Sympc0{\operatorname{Symp^c_0}}

\def\Ham{\operatorname{Ham}}

\def\length{\operatorname{length}}

\def\ind{\operatorname{ind}}

\def\can{\operatorname{can}}

\def\top{\operatorname{top}}
\def\fibre{\operatorname{fibre}}
\def\sphere{\operatorname{sphere}}
\def\HF{\operatorname{HF}}
\def\CF{\operatorname{CF}}
\def\H{\operatorname{H}}
\def\dia{\operatorname{diam}}
\def\const{\operatorname{const}}
\def\can{\operatorname{can}}

\def\gb{\beta}
\def\gg{\gamma}
\def\gd{\delta}
\def\eps{\epsilon}

\def\gf{\varphi}

\def\gl{\lambda}
\def\go{\omega}
\def\gs{\sigma}
\def\gt{\vartheta}

\def\ca{{\mathcal A}}

\def\ce{{\mathcal E}}

\def\cj{{\mathcal J}}
\def\cl{{\mathcal L}}
\def\cm{{\mathcal M}}

\def\cp{{\mathcal P}}

\def\bJ{{\mathbf J}}

\def\NN{\mathbbm{N}}

\def\QQ{\mathbbm{Q}}
\def\RR{\mathbbm{R}}

\def\ZZ{\mathbbm{Z}}

\def\RP{\operatorname{\mathbbm{R}P}}  
\def\CP{\operatorname{\mathbbm{C}P}}
\def\HP{\operatorname{\mathbbm{H}P}}
\def\Ca{\operatorname{\mathbbm{C}\mathbbm{a}P}^2}

\def\SO{\operatorname{SO}}
\def\SU{\operatorname{SU}}

\def\pp{\partial}

\def\ra{\rightarrow}

\def\ni{\noindent}
\def\b{\bigskip}
\def\m{\medskip}

\def\id{\mbox{id}}

\def\proof{\noindent {\it Proof. \;}}
\begin{document}

\title{Fiberwise volume growth via Lagrangian intersections}

\author{Urs Frauenfelder}
\address{(U.\ Frauenfelder) Department of Mathematics, Hokkaido University,
Sapporo 060-0810, Japan}
\email{urs@math.sci.hokudai.ac.jp}
\author{Felix Schlenk}
\address{(F.\ Schlenk) Mathematisches Institut,
Universit\"at Leipzig, 04109 Leipzig, Germany}
\email{schlenk@math.uni-leipzig.de}
\thanks{This work was supported by the 
Japanese Society for Promotion of Science (UF), 
and the Deutsche Forschungsgemeinschaft (FS)}

\date{\today}

\begin{abstract}
We consider Hamiltonian diffeomorphisms $\gf$ of the unit
cotangent bundle over a closed Riemannian manifold 
$(M,g)$ which extend to Hamiltonian diffeomorphisms of $T^*M$
equal to the time-$1$-map of the geodesic flow for $|p| \ge 1$.
For such diffeomorphisms we
establish uniform lower bounds for the fiberwise volume
growth of $\gf$ which were previously known for geodesic flows
and which depend only on $(M,g)$ or on the homotopy type of $M$.
More precisely, we show that for each $q \in M$ the volume growth
of the unit ball in $T_q^*M$ under the iterates of $\gf$ is 
at least linear if $M$ is rationally elliptic, 
is exponential if $M$ is rationally hyperbolic,
and is bounded from below by the growth of the fundamental group
of $M$.
In the case that all geodesics of $g$ are closed, 
we conclude that the slow volume growth of every symplectomorphism in the
symplectic isotopy class of the Dehn--Seidel twist is at least
$1$, completing the main result of~\cite{FS:GAFA}.
The proofs use the Lagrangian Floer homology of $T^*M$ and the
Abbondandolo--Schwarz isomorphism from this homology to the
homology of the based loop space of $M$.
\end{abstract}

\maketitle

\markboth{{\rm Fiberwise volume growth via Lagrangian intersections}}{{}} 


\section{Introduction and main results}  

\ni
\subsection{Topological entropy and volume growth}
The topological entropy $h_{\top} (\gf)$ of a compactly supported
$C^1$-diffeomorphism $\gf$ of a smooth manifold $X$ is a basic numerical
invariant measuring the orbit structure complexity of $\gf$.
There are various ways of defining $h_{\text{top}} (\gf)$, see \cite{KH}.
If $\gf$ is $C^\infty$-smooth, 
a geometric way was found by Yomdin and Newhouse in their seminal works
\cite{Y} and \cite{N}:
Fix a Riemannian metric $g$ on $X$.  
For $j \in \left\{ 1, \dots, \dim X \right\}$
denote by $\Sigma_j$ the set of smooth compact 
(not necessarily closed) $j$-dimensional submanifolds of $X$, 
and by $\mu_g(\gs)$ the volume of $\gs \in \Sigma_j$
computed with respect to the measure on $\gs$ induced by $g$. 
The {\it $j$'th volume growth}\, of $\gf$ is defined as 
\[
v_j (\gf) \,=\, 
 \sup_{\gs \in \Sigma_j} \liminf_{m \ra \infty} \frac 1m
    \log \mu_g \left( \gf^m ( \gs ) \right) ,
\]
and the {\it volume growth}\, of $\gf$ is defined as 
\[                   
v (\gf) \,=\, \max_{1 \le j \le \dim X} v_j(\gf) .
\]
Newhouse proved in \cite{N} that $h_{\text{top}} (\gf) \le v (\gf)$,
and Yomdin proved in \cite{Y} that $h_{\top} (\gf) \ge v (\gf)$ 
provided that $\gf$ is $C^\infty$-smooth, so that
\begin{equation}  \label{e:yomdin}
h_{\top} (\gf) \,=\, v (\gf) \quad \text{if $\gf$ is $C^\infty$-smooth}.
\end{equation}
The topological entropy measures the {\it exponential}\, growth
rate of the orbit complexity of a diffeomorphism.
It therefore vanishes for many interesting dynamical systems.
Following \cite{KT,FS:GAFA} we thus also consider the {\it
$j$'th slow volume growth} 
\[
s_j (\gf) \,=\, 
 \sup_{\gs \in \Sigma_j} \liminf_{m \ra \infty} \frac {1}{\log m} 
    \log \mu_g \left( \gf^m ( \gs ) \right) 
\]
and the {\it slow volume growth}
\[
s (\gf) \,=\, \max_{1 \le j \le \dim X} s_j(\gf) .
\]
It measures the {\it polynomial}\, volume growth of the iterates of
the most distorted smooth $j$-dimensional family of initial
data.
Note that $v_j (\gf)$, $v (\gf)$, $s_j(\gf)$, $s(\gf)$ do not
depend on the choice of $g$, and that $v_{\dim X} (\gf) = s_{\dim X} (\gf) =0$.

\m
The aim of this paper is to give uniform lower estimates of
localized versions of 
$v(\gf)$ and $s(\gf)$ for certain symplectomorphisms of cotangent bundles.
We consider a smooth closed $d$-dimensional Riemannian manifold
$(M,g)$ and the cotangent bundle $T^*M$ over $M$ endowed with
the induced Riemannian metric $g^*$ and the 
standard symplectic form $\go = \sum_{j=1}^d dp_j \wedge dq_j$.
We abbreviate
\[
D(r) \,=\, \left\{ (q,p) \in T^*M \mid |p| \le r \right\}
\quad \text{ and } \quad
D_q(r) \,=\, T_q^*M \cap D(r) .
\]
Let $\gf$ be a $C^1$-smooth symplectomorphism of $\left( T^*M,
\go \right)$ which preserves $D(r)$.
If $\gf$ is $C^\infty$-smooth, \eqref{e:yomdin} says that the
maximal orbit complexity of $\gf |_{D(r)}$ is already
contained in the orbit of a single submanifold of $D(r)$.
Usually, lower estimates of the topological entropy do not give any
information on the dimension or the location of such a submanifold.
An attempt to localize such submanifolds for symplectomorphisms
was made in \cite{FS:GAFA}, where we considered 
{\it Lagrangian}\, submanifolds only.
In this paper we further localize and consider for $\gf$ as
above the 
{\it fiberwise volume growth}
\[
\hat{v}_{\fibre} (\gf;r) \,=\, 
 \sup_{q \in M} \liminf_{m \ra \infty} \frac 1m
    \log \mu_{g^*} \bigl( \gf^m \left( D_q(r) \right) \bigr)
\]
and the {\it slow fiberwise volume growth}
\[
\hat{s}_{\fibre} (\gf;r) \,=\, 
 \sup_{q \in M} \liminf_{m \ra \infty} \frac{1}{\log m}
    \log \mu_{g^*} \bigl( \gf^m \left( D_q(r) \right) \bigr) .
\]
In fact, we shall give uniform lower estimates of the (slow) volume
growth of {\it each}\, fibre by considering the 
{\it uniform fiberwise volume growth}
\[
\check{v}_{\fibre} (\gf;r) \,=\, 
 \inf_{q \in M} \liminf_{m \ra \infty} \frac 1m
    \log \mu_{g^*} \bigl( \gf^m \left( D_q(r) \right) \bigr)
\]
and the {\it uniform slow fiberwise volume growth}
\[
\check{s}_{\fibre} (\gf;r) \,=\, 
 \inf_{q \in M} \liminf_{m \ra \infty} \frac{1}{\log m}
    \log \mu_{g^*} \bigl( \gf^m \left( D_q(r) \right) \bigr) .
\]
Writing $v(\gf;r) = v \left( \gf |_{D(r)} \right)$ and so on,
we clearly have
\begin{eqnarray}  
v(\gf;r) &\ge& v_d(\gf;r) \,\ge\, \hat{v}_{\fibre} (\gf;r)
\,\ge\, \check{v}_{\fibre} (\gf;r) , \label{ines:v} \\ [.5em]
s(\gf;r) &\ge& s_d(\gf;r) \,\ge\, \hat{s}_{\fibre} (\gf;r)
\,\ge\, \check{s}_{\fibre} (\gf;r) . \label{ines:s}
\end{eqnarray}
We shall obtain lower estimates of $\check{v}_{\fibre} (\gf;r)$ and
$\check{s}_{\fibre} (\gf;r)$ in terms of the growth of certain
homotopy-type invariants of $M$.
We introduce these concepts of growth right now.

\subsection{Rationally elliptic and hyperbolic manifolds}  \label{ss:rational}
Let $M$ be a closed connected manifold whose fundamental group $\pi_1(M)$ is {\it finite}.
Such a manifold $M$ is said to be {\it rationally elliptic}\, if
the total rational homotopy $\pi_* (M) \otimes \QQ$ is finite
dimensional, and $M$ is said to be  
{\it rationally hyperbolic}\, if the integers
\begin{equation}  \label{def:rho}
\rho_m (M) \,:=\, \sum_{j=0}^m \dim \pi_j (M) \otimes \QQ
\end{equation}
grow exponentially, i.e., there exists $C>1$ such that 
$\rho_m(M) \ge C^m$ for all large enough $m$.
It is shown in \cite{FH,FHT} that every closed manifold $M$ with $\pi_1(M)$ finite 
is either rationally elliptic or rationally hyperbolic. The
proof is based on Sullivan's minimal models. 

\b
\ni
{\bf Examples.}
We assume that $M$ is simply connected. In dimensions $2$ and
$3$, the standard sphere is the only such manifold up to diffeomorphism in
view of the proof of the Poincar\'e conjecture; 
it is rationally elliptic.
In higher dimensions, ``most'' simply connected mani\-folds are rationally hyperbolic.
In dimension $4$, the simply connected rationally elliptic
manifolds up to homeomorphism are 
\[
S^4, \quad
\CP^2, \quad
S^2 \times S^2, \quad
\CP^2 \# \overline{\CP}^2, \quad
\CP^2 \# \CP^2 ,
\]
see \cite[Lemma~5.4]{Pa}, 
and in dimension $5$ the simply connected rationally elliptic
manifolds up to diffeomorphism are 
\[
S^5, \quad
S^2 \times S^3, \quad
S^2 \ltimes S^3, \quad 
\SU (3) / \SO (3), 
\]
where $S^2 \ltimes S^3$ is the nontrivial $S^3$-bundle over $S^2$,
see \cite{PP}.
\diam

We refer to \cite{FH,FHT,Pa,PP} and the references therein for more
information on rationally elliptic and hyperbolic manifolds.

\subsection{Growth of infinite fundamental groups}  \label{s:gr}
Assume now that $\pi_1(M)$ is {\it infinite}.
Since $M$ is closed, $\pi_1(M)$ is then an infinite finitely presented group.
Consider, more generally, an infinite finitely generated group $\Gamma$. 
The growth function $\gamma_S$ associated with 
a finite set $S$ of generators of $\Gamma$ is defined as follows: 
For each positive integer $m$, let $\gamma_S (m)$ 
be the number of distinct group elements which can be written as
words with at most $m$ letters from $S \cup S^{-1}$.
As is easy to see, \cite{Mi}, the limit
\[
\nu (S) \,:=\, \lim_{m \to \infty} \frac{\log \gg_S (m)}{m}
\,\in\, [0,\infty)
\]
exists, and the property $\nu (S) > 0$ is independent of the choice of $S$.
In this case, 
the group $\Gamma$ is said to have {\it exponential growth}.

\begin{example*}
{\rm
If a closed manifold $M$ admits a Riemannian metric of negative
sectional curvature,
then $\pi_1(M)$ has exponential growth, \cite{Mi}.
For $d=2$, the converse to this statement holds true, while for
$d \ge 3$ there are closed $d$-dimensional manifolds for which
$\pi_1(M)$ has exponential growth and which carry no
Riemannian metric with negative sectional curvature, 
see~\cite[Corollaire~III.10]{Ma} and \cite[p.~190]{Harpe}.
\diam
}
\end{example*}

If $\Gamma$ does not have exponential growth, $\Gamma$ is said to have {\it subexponential growth}.
In this case, the degree of polynomial growth 
\[
s (\Gamma) \,:=\, \limsup_{m \to \infty} \frac{\log \gg_S (m)}{\log m}
\,\in\, [0,\infty] 
\]
is independent of the choice of $S$.
The group $\Gamma$ has {\it intermediate growth}\, 
if $s (\Gamma) = \infty$, and {\it polynomial growth}\, 
if $s (\Gamma) < \infty$.
While there are finitely generated groups of intermediate growth, \cite{Gri,Harpe},
it is still unknown whether there are finitely presented groups of
intermediate growth;
we shall thus not consider a more refined notion of growth for such
groups.
According to a theorem of Gromov, \cite{Gr81}, the group $\Gamma$ has
polynomial growth if and only if 
$\Gamma$ has a nilpotent subgroup $\Delta$ of finite index.
Let $\left( \Delta_k \right)_{k \ge 1}$ be its lower central series
defined inductively by
$\Delta_1 = \Delta$ and $\Delta_{k+1} = [ \Delta, \Delta_k]$. 
Then
\[
s (\Gamma ) \,=\, 
\sum_{k \ge 1} k \dim \bigl( \left(\Delta_k / \Delta_{k+1} \right) \otimes_\ZZ \QQ \bigr) ,
\]
see \cite{Wolf,Gui,Bass,Tits}.
We in particular see that $s(\Gamma)$ is a positive integer.

\begin{examples*}
{\rm
(i) 
For the fundamental group of the torus, $s (\ZZ^d) = 1 \cdot d = d$.

\s
(ii) 
For the Heisenberg group
\[
\Gamma \,=\, 
\left\{
         \begin{pmatrix}
               1 & 0  & 0  &    \\
               x & 1  & 0  &    \\
               z & y  & 1  &     
         \end{pmatrix}      
\Bigg| \;
x, y, z \in \ZZ
\right\}
\]
we have $\Delta_1 = \Gamma$ and  
$\Delta_2 = \left\{ M(x,y,z) \in \Gamma \mid x=y=0 \right\} \cong
\ZZ$ and $\Delta_k = \{ e \}$ for $k \ge 3$, so that $s (\Gamma) = 1 \cdot
2 + 2 \cdot 1 =4$.
}
\diam
\end{examples*}
Much more information on growth of finitely generated groups 
can be found in \cite{Harpe}.

\s
Consider now a closed connected Riemannian manifold $(M,g)$ such
that $\pi_1(M)$ has exponential growth.
Given a finite set $S$ of generators 
let $l (S,g)$ be the smallest real number such that for each 
$q \in M$ 
each generator $s \in S$ of $\pi_1(M) \cong \pi_1(M,q)$ can be
represented by a smooth loop based at $q$ of length no more than
$l(S,g)$.
Then set
\[
\nu (M,g) \,=\, \sup \frac{\nu(S)}{l(S,g)}
\]
where the supremum is taken over all finite sets $S$ generating
$\pi_1(M)$.

\subsection{Main result}
Consider again a closed Riemannian manifold $(M,g)$, and let 
$H \colon [0,1] \times T^*M \ra \RR$ be a $C^2$-smooth 
Hamiltonian function meeting the following assumption:
There exists $r_H>0$ and a function
$f \colon [0,\infty) \ra \RR$ with $f'(r_H) \neq 0$ such that 
\begin{equation}  \label{ge:r}
H(t,q,p) \,=\, f \left( |p| \right) \quad \text{ for }\,|p| \ge r_H. 
\end{equation}
The Hamiltonian flow $\gf_H^t$ of the time-dependent vector
field $X_H$ given by  
$\go \left( X_{H_t}, \cdot \right) = -dH_t \left( \cdot \right)$
is defined for all $t \in [0,1]$.
We abbreviate $\gf_H = \gf_H^1$.

\b
\ni
{\bf Theorem~1.}
{\it
Consider a closed Riemannian manifold $(M,g)$,
and let $H \colon [0,1] \times T^*M \ra \RR$ be a $C^2$-smooth Hamiltonian 
function satisfying~\eqref{ge:r}.
\begin{itemize}
\item[(i)]
Assume that $\pi_1(M)$ is finite.
If $M$ is rationally elliptic, 
then 
\[
\check{s}_{\fibre} (\gf_H;r_H) \,\ge\, 1 .
\]
If $M$ is rationally hyperbolic,
then 
\[
\check{v}_{\fibre} (\gf_H;r_H) \,\ge\, f'(r_H) \, r_H \, C_1
\]
for some positive constant $C_1$ depending only on $(M,g)$.

\item[(ii)]
Assume that $\pi_1(M)$ is infinite.
If $\pi_1(M)$ has subexponential growth, 
then 
\[
\check{s}_{\fibre} (\gf_H;r_H) \,\ge\, s \left( \pi_1(M) \right) .
\]
If $\pi_1(M)$ has exponential growth, 
then 
\[
\check{v}_{\fibre} (\gf_H;r_H) \,\ge\, 2 \, f'(r_H) \, r_H \, \nu (M,g) .
\]
\end{itemize}
}

\b
\ni
{\bf Discussion~1.}
{\rm
(i) 
There are rationally elliptic manifolds for which all the
numbers in \eqref{ines:s} are~$1$, see Discussion~2\:(i)
below.
The first statement in (i) is thus sharp. 
For the flat torus $T^d$ and $H = \frac 12 p^2$, we have
$\check{s}_{\fibre} (\gf_H;r) = s (\ZZ^d) =d$ for all $r>0$, so that the
first statement in (ii) is sharp.

\s
(ii) 
If $H$ is $C^\infty$-smooth, then 
$h_{\top} (\gf_H;r_H) \ge \check{v}_{\fibre}(\gf_H;r_H)$ 
by Yomdin's theorem and \eqref{ines:v}, so that the second
statements in (i) and (ii) yield positive lower bounds for  
$h_{\top} (\gf_H;r_H)$.
These bounds imply and are implied by the estimates
\[
h_{\top} (g) \,\ge\, C_1(M,g) 
\quad \text{ and } \quad
h_{\top} (g) \,\ge\, 2 \,\nu (M,g) 
\]
for the topological entropy of the geodesic flow on the unit
sphere bundle $\pp D(1)$ of a $C^\infty$-smooth Riemannian
metric $g$
on a rationally hyperbolic manifold or a manifold whose
fundamental group has exponential growth.
The first of these estimates was found by 
Gromov and Paternain (see \cite[Corollary~5.21]{Pa}),
and the second estimate is a version of Dinaburg's theorem,
which holds for $C^2$-smooth $g$, 
(see \cite{Dina} or \cite[Theorem~5.18]{Pa}).

\s
(iii) 
Theorem~1 extends well-known results from the study of geodesic
flows, see \cite[Corollary~3.9 and Chapter~5]{Pa}:
These results imply Theorem~1 if there exists an $\eps >0$ such
that $H = \frac 12 |p|^2$ on $D(r_H) \setminus D(r_H-\eps)$.
In this situation, these results as well as Theorem~1 itself
imply the second statements in (i) and (ii) with 
$\check{v}_{\fibre} (\gf_H;r_H)$ replaced by the {\it uniform
spherical volume growth}\,
\[
\check{v}_{\sphere} (\gf_H;r_H) \,=\,
 \inf_{q \in M} \liminf_{m \ra \infty} \frac 1m
    \log \mu_{g^*} \bigl( \gf^m_H \left( \pp D_q(r_H) \right)
\bigr) .
\]

\s
(iv)  
As the identity mapping illustrates, the assumption $f'(r) \neq
0$ in \eqref{ge:r} is essential.
\diam
}

\subsection{Volume growth in the component of the
Dehn--Seidel twist}
A {\it $P$-manifold}\, is a connected Riemannian manifold all of
whose geodesics are periodic. 
Such manifolds are closed, and as we shall see in Section~\ref{s:P}, 
every \text{$P$-manifold} different from $S^1$ is rationally elliptic.
The known \text{$P$-manifolds} are the compact rank one symmetric spaces
\[
S^d, \quad \RP^d, \quad \CP^n, \quad \HP^n, \quad \Ca 
\]       
with their canonical Riemannian structures, their Riemannian
quotients (which are all known), and so-called Zoll manifolds,
which are modelled on spheres.
It is an open problem whether there are other $P$-manifolds. 
More information on $P$-manifolds can be found in
\cite{B,FS:GAFA}, Section~10.10 of \cite{Be} and
Section~\ref{s:P} below.

Let $(M,g)$ be a $P$-manifold.
It is known that the unit-speed geodesics of $(M,g)$ admit a common
period, and we shall assume $g$ to be scaled such that the minimal
common period is $1$.  
We choose a smooth function $f \colon [0,\infty) \ra [0,\infty)$
such that 
\begin{equation}  \label{con:f}
f (r) = \tfrac 12 r^2 \;\text{ near } 0
\qquad \text{and} \qquad
f'(r) =1  \;\text{ for } r \ge 1 ,
\end{equation}
and following \cite{A,S1,S2} 
we define the (left-handed) Dehn--Seidel twist $\gt_f$ to
be the time-1-map of the Hamiltonian flow generated by 
$f \left( |p|\right)$.
Since $(M,g)$ is a $P$-manifold, $\gt_f$ is the identity on
$T^*M \setminus T^*_1M$, so that $\gt_f$ is a compactly
supported symplectomorphism, $\gt_f \in \Symp^c \left( T^*M
\right)$.
We shall write $\gt$ for any map $\gt_f$ with $f$ satisfying \eqref{con:f}.
The class $\left[ \gt \right]$ in the symplectic mapping
class group $\pi_0 \left( \Symp^c \left( T^*M \right) \right)$
clearly does not depend on $f$.
Given $\gf \in \Symp^c \left( T^*M \right)$ we set 
$\check{s}_{\fibre} (\gf) = \check{s}_{\fibre} (\gf;r)$
where $r >0$ is any number such that $\gf$ is supported in $D(r)$. 

\b
\ni
{\bf Corollary~1.}
{\it
Let $(M,g)$ be a $P$-manifold, and let $\gt$ be a twist on
$T^*M$.
If $\gf \in \Symp^c \left( T^*M \right)$ is such that $[\gf]=
\left[ \gt^m \right] \in \pi_0 \left( \Symp^c \left( T^*M
\right) \right)$ for some $m \in \ZZ \setminus \{ 0 \}$, then
$\check{s}_{\fibre} (\gf) \ge 1$.
}

\b
\ni
{\bf Discussion~2.}
{\rm
(i)
A computation given in \cite[Proposition~2.2\:(i)]{FS:GAFA} shows that 
$s \left( \gt_f^m \right) = \check{s}_{\fibre} \left( \gt_f^m \right) = 1$ 
for every $m \in \ZZ \setminus \{0\}$ and every $f$.
Corollary~1 is thus sharp and shows that twists
minimize slow volume growths in their symplectic isotopy class.

\s
(ii)
Assume that $\gf \in \Symp^c \left( T^*M \right) \setminus \{
\id \}$
is such that $[\gf] = [\gt^0] = [ \id ]$.
If $d \ge 2$, this means that $\gf$ is a non-identical compactly supported
Hamiltonian diffeomorphism of $T^*M$.
If the support of $\gf$ misses some fibre, then $\check{s}_{\fibre} (\gf) = 0$.
On the other hand, 
combining a result in \cite{FS} with the arguments in \cite{P1}
one finds that $s_1(\gf) \ge 1$, see \cite{FS2}.
It is not hard to construct examples with $s(\gf) = s_j(\gf) =
1$ for all $j \in \{ 1, \dots, 2d-1 \}$, see \cite{FS2}.   

\s
(iii)
It was proved in \cite[Corollary~4.5]{S2} that the class
$[\gt]$ of a twist generates an infinite cyclic subgroup of
$\pi_0 \left( \Symp^c \left( T^*M \right) \right)$.
Corollary~1 yields another proof of this.
            
\s
(iv)
Corollary~1 was proved in \cite{FS:GAFA} for all known (but
possibly not for all) $P$-manifolds.
The proof there only used Lagrangian Floer homology and a
symmetry argument for $\left( \CP^n, g_{\can} \right)$.
}

\b
\ni
{\bf Acknowledgements.}
We wish to thank the referee of our paper \cite{FS:GAFA} 
for suggesting the present approach to Corollary~1.
Another important stimulus was Gabriel Paternain's beautiful
lecture on topological entropy at Max Planck Institute Leipzig
in November~2004.
Most of this paper was written at the turn of the year 2004/2005 at Hokkaido 
University at Sapporo.
We cordially thank Kaoru Ono for valuable discussions and for
his generous hospitality, 
and Theo B\"uhler, Leonardo Macarini and Matthias Schwarz for helpful comments.

\section{Proofs}

\subsection{Proof of Theorem~1}
Our proof is along the following lines.
Using an idea from~\cite{FS:GAFA} we first show that fiberwise
volume growth is a consequence of the growth of the dimension of
certain Floer homology groups.
Applying the isotopy invariance of Floer homology and
a recent result of Abbondandolo and Schwarz, these
homology groups are seen to be isomorphic to the homology of the
space of based loops in $M$ not exceeding a certain length.
Their dimension can be estimated from below by results of Gromov and
Serre if $\pi_1(M)$ is finite and by elementary considerations
if $\pi_1(M)$ is infinite.

\b
\ni
{\bf Step~1. Volume growth via Floer homology} 

\m
\ni
{\bf Geometric set-up.}
Let $(M,g)$ be a closed Riemannian manifold, and let $H \colon
[0,1] \times T^*M \to \RR$ be as in Theorem~1.
We abbreviate $\gb = f'(r_H)$, and
we assume without loss of generality that $\gb >0$.
Fix $\eps >0$ and let
$K \colon [0,1] \times T^*M \to \RR$ be a $C^\infty$-smooth
function such that
\begin{equation}  \label{e:K}
K(t,q,p) = \gb |p|^2 \quad \text{ if }\, 
|p| \ge r_H + \eps .
\end{equation}
For $m =1, 2, \dots$ we recursively define 
$K_t^m=K_t+K_t^{m-1} \circ \gf_K^t$.
Then $\gf_K^m=\gf_{K^m}$.
For $q_0, q_1 \in M$ let $\Omega^1 \left( T^*M,q_0,q_1\right)$
be the space of all paths $x \colon [0,1] \to T^*M$ of Sobolev
class $W^{1,2}$ such that $x(0) \in T_{q_0}^*M$ and $x(1) \in
T_{q_1}^*M$.
This space has a canonical Hilbert manifold structure, \cite{Kli}.
The action functional of classical mechanics 
$\ca_{K^m} \colon \Omega^1 \left( T^*M,q_0,q_1 \right) \to \RR$
associated with $K^m$ is defined as
\[
\ca_{K^m}(x) \,=\, 
\int_0^1 \bigl( \gl \left(\dot{x}\right) - K^m (t,x) \bigr) dt,
\]
where $\gl = \sum_{j=1}^d p_j \,dq_j$ is the canonical $1$-form on $T^*M$.
This functional is $C^\infty$-smooth, and its
critical points are precisely the elements of
the space $\cp \left( q_0,q_1,K^m \right)$ of $C^\infty$-smooth
paths $x \colon [0,1] \to T^*M$ solving
\[
\dot{x}(t) = X_{K^m}\left( x(t) \right), \,\,t \in [0,1], \quad \;
x(j) \in T_{q_j}^*M, \,\,j=0,1.
\]
Notice that the elements of $\cp \left( q_0,q_1,K^m \right)$
correspond to the intersection points of 
$\gf_K^m \left( T^*_{q_0}M \right)$ and $T^*_{q_1}M$ 
via the evaluation map $x \mapsto x(1)$.

Fix now $q_0 \in M$, and let $V (q_0, K^m)$ be the
set of those $q_1 \in M$ for which 
$\gf_K^m \left( D_{q_0} (r_H + 2 \eps) \right)$ and
$D_{q_1} (r_H + 2 \eps)$ intersect transversely. 

\begin{lemma}\label{l:full}
The set $V (q_0, K^m)$ is open and of full measure in $M$.
\end{lemma}

\proof
Since $D_{q_0} (r_H+2\eps)$ is compact, $V (q_0, K^m)$ is open.
Applying Sard's Theorem to the projection 
$\gf_K^m \left( D_{q_0} (r_H+2\eps) \right) \to M$ 
one sees that $V (q_0, K^m)$ has full measure in $M$.
\proofend

Consider the ``annuli-bundle''
\[
A(\eps) \,=\, \left\{ (q,p) \in T^*M \mid r_H + \eps \le |p| \le
r_H+2\eps \right\} ,
\]
set $A_{q_0}(\eps) = A(\eps) \cap T_{q_0}^*M$, and let 
$W (q_0,K^m)$ be the set of those 
$q_1 \in V (q_0, K^m)$ for which
$\gf_K^m \left( A_{q_0} (\eps) \right) \cap A_{q_1} (\eps)$ is
empty.
For $q_1 \in W (q_0,K^m)$ the set
\[
\cp \left( q_0,q_1,K^m, r_H + 2\eps \right)
\,:=\, \left\{ x \in \cp \left( q_0,q_1,K^m \right)
\mid x \subset D \left( r_H + 2 \eps \right) \right\}
\] 
is finite and contained in $D \left( r_H + \eps \right)$.
Since the first Chern class of $\left( T^*M, \go \right)$
vanishes,
each $x \in \cp \left( q_0,q_1,K^m, r_H + 2\eps \right)$ comes
with an integral index, which in case of a geodesic
Hamiltonian agrees with the Morse index of the corresponding
geodesic path, see \cite{We,AS}.

\b
\ni
{\bf Floer homology.}
Floer homology for Lagrangian intersections was invented by Floer
in a series of seminal papers, \cite{F2,F3,F1,F4}.
We shall use a version of Floer homology described in \cite{KS,FS:GAFA}.
In the above situation, we define the $k^{th}$ Floer chain group
$\CF_k \left( q_0,q_1,K^m,r_H+2\eps \right)$ 
as the finite-dimensional $\QQ$-vector space freely
generated by the elements of $\cp \left( q_0,q_1,K^m, r_H +
2\eps \right)$ of index $k$, and the full Floer chain group
as
\[
\CF_* \left( q_0,q_1,K^m,r_H+2\eps \right) 
\,=\,
\bigoplus_{k\in \ZZ} \CF_k \left( q_0,q_1,K^m,r_H+2\eps \right) .
\]
In order to define the Floer boundary operator, we follow 
\cite{CFH,V2,BPS} and consider the set $\cj$ of $t$-dependent 
smooth families $\bJ = \{ J_t\}$, $t \in [0,1]$, 
of $\go$-compatible almost complex structures on $D \left(r_H+2\eps\right)$
such that $J_t$ is convex and independent of $t$ on $A(\eps)$.
This in particular means that $J$ is invariant under the local
flow of the Liouville vector field 
$Y = \sum p_j \pp p_j$ on $A(\eps)$.
For $\bJ \in \cj$, for smooth maps $u$ from the strip $S = \RR
\times [0,1]$ to $D \left( r_H + 2 \eps \right)$, and for $x^\pm
\in \cp \left( q_0,q_1,K^m,r_H+2\eps \right)$ 
consider Floer's equation
\begin{eqnarray}  \label{floer}
 \left\{       
  \begin{array}{lcr}
    \pp_s u+J_t(u) \left( \pp_t u - X_{K_t^m}(u) \right) = 0,  \\ [0.2em]
    u(s,j) \in T^*_{q_j}M,\,\, j=0,1 , \\ [0.2em]
    \displaystyle{\lim_{s \to \pm \infty}} u(s,t) \,=\, x^\pm (t) 
             \;\text{ uniformly in } t .
  \end{array}
 \right.
\end{eqnarray}

\begin{lemma}  \label{l:convex}
Solutions of \eqref{floer} are contained in $D (r_H + \eps)$.
\end{lemma}

\ni
{\it Sketch of proof.}
Since $q_1 \in W(q_0,K^m)$,
\begin{equation}  \label{e:limit}
\lim_{s \to \pm \infty} u(s,t) \,=\, x^\pm (t) \subset D(r_H + \eps) .
\end{equation}
In view of the strong maximum principle, the lemma follows from the
convexity of $J$ on $A(\eps)$ and from \eqref{e:limit} 
together with the fact that the special form
\eqref{e:K} of $K$ on $A(\eps)$ implies 
$\go \left( Y, J X_{K^m} \right) =0$, 
cf.\ \cite{KS,FS:GAFA}.
\proofend

We denote the set of solutions of~\eqref{floer} by 
$\cm \left( x^-,x^+,K^m,r_H+2\eps\right)$.
Lemma~\ref{l:convex} is an important ingredient to establish the
compactness of this set.
The other ingredient is that there is no bubbling-off of 
$\bJ$-holomorphic spheres or discs.
Indeed, $[\go]$ vanishes on $\pi_2 \left( T^*M \right)$ because
$\go = d \gl$ is exact, and $[\go]$ vanishes on 
$\pi_2 \bigl( T^*M , T_{q_j}^*M \bigr)$ because $\gl$ vanishes on
$T_{q_j}^*M$, $j=0,1$.

It is shown in \cite{AS} that $\cm \left( x^-,x^+,K^m,r_H+2\eps\right)$
admits a coherent orientation.
For a generic choice of $\bJ \in \cj$ the Floer boundary operators
\[
\pp_k \left( \bJ \right) \colon
\CF_k \left( q_0,q_1,K^m,r_H+2\eps \right) \,\to\,
\CF_{k-1} \left( q_0,q_1,K^m,r_H+2\eps \right)
\]
can now be defined in the usual way. 
The Floer homology groups with rational coefficients
\[
\HF_k \left( q_0,q_1,K^m,r_H+2\eps \right)
\,:=\, \ker \pp_k (\bJ) / \im \pp_{k+1} (\bJ)
\] 
do not depend on the choice of $\bJ$ up to natural isomorphisms,
nor do they alter if we add to $K^m$ a function supported in
$[0,1] \times D (r_H+\eps)$.
The function $G_\gb^m := m \gb |p|^2$ is such a function, so
that
\begin{equation}  \label{e:isotopy}
\HF_* \left( q_0,q_1,K^m,r_H+2\eps \right)
\,\cong \, 
\HF_* \left( q_0,q_1,G_\gb^m,r_H+2\eps \right)
\end{equation}
provided that $q_1$ also belongs to $V (q_0, G_\gb^m)$.
\begin{remark}  
{\rm
For most applications of Floer homology found so far,
it suffices to work with the coefficient field $\ZZ_2$. In order
to prove the second statement in Theorem~1\;(i) it will be
important that we can work with rational coefficients,
see Remark~\ref{r:rational} below.
}
\end{remark}

\begin{proposition}  \label{p:floer}
Assume that $q_1 \in W (q_0,G_\gb^m)$.
\begin{itemize}
\item[(i)]
Assume that $M$ is simply connected.
If $M$ is rationally elliptic, 
then there exists a constant $c_1 >0$ such that
\[
\dim \HF_* \left( q_0,q_1, G_\gb^m, r_H+2\eps \right)
\,\ge\, \left( c_1 \gb r_H \right) m   
\quad \; \text{ for all large enough }  m \in \NN .
\]
If $M$ is rationally hyperbolic, 
then there exists a constant $C_1 >0$ 
depending only on $(M,g)$ such that
\[
\dim \HF_* \left( q_0,q_1, G_\gb^m, r_H+2\eps \right)
\,\ge\, e^{C_1 \gb r_H m}  
\quad\; \text{ for all large enough }  m \in \NN . 
\] 
\item[(ii)]
Assume that $\pi_1(M)$ is infinite. 
Fix a set $S$ of generators, let $\gg_S$ be the
corresponding growth function,
and let $l(S,g)$ be as defined in Section~\ref{s:gr}.
Then
\[
\dim \HF_* \left( q_0,q_1, G_\gb^m, r_H+2\eps \right)
\,\ge\, \gg_S \left( \lfloor (2 \gb r_H) m / l(S,g) \rfloor \right) 
\quad\; \text{ for all large enough }  m \in \NN . 
\] 
\end{itemize}
\end{proposition}

Here, $\lfloor r \rfloor = \max \left\{ n \in \ZZ  \mid n \le r \right\}$.
Proposition~\ref{p:floer} will be proved in the next two steps.
In the remainder of this step we show

\b
\ni
{\bf Proposition~\ref{p:floer} $\Longrightarrow$ Theorem~1.}
We show this implication for rationally elliptic manifolds.
The other implications are shown in a similar way.

We assume first that $M$ is a rationally elliptic manifold which
is simply connected.
Let $K$ be as before, and pick $q_1 \in W (q_0, K^m) \cap
W(q_0,G_\gb^m)$.
Since the generators of $\CF_* \left( q_0,q_1, K^m, r_H+2\eps \right)$
correspond to 
$\gf_K^m \left( D_{q_0} \left( r_H+\eps \right) \right) 
\cap D_{q_1} \left( r_H+\eps\right)$, we find together
with \eqref{e:isotopy} and Proposition~\ref{p:floer}\,(i) that
\begin{eqnarray}
\# \bigl( \gf_K^m \left( D_{q_0} (r_H+\eps )
\right) \cap D_{q_1} (r_H+\eps ) \bigr) \notag
&=&
\dim \CF_* \left( q_0,q_1, K^m,r_H+2\eps \right) \notag \\
&\ge& 
\dim \HF_* \left( q_0,q_1, K^m,r_H+2\eps \right) \notag \\
&=& 
\dim \HF_* \left( q_0,q_1, G_\gb^m,r_H+2\eps \right) \notag \\
&\ge& 
\left( c_1 \gb r_H \right) m \label{e:nC}
\end{eqnarray}
for $m$ large enough.
Recall now that $\gb = f'(r_H)$. We thus find a sequence
$\eps_i \to 0$ and a sequence $K_i \colon [0,1] \times T^*M
\to \RR$ of $C^\infty$-smooth functions such that 
\[
K_i(t,q,p) = \gb |p|^2 \quad \text{ if }\, 
|p| \ge r_H + \eps_i 
\]
and such that
\begin{itemize}
\item[(K1)]
$K_i |_{D(r_H+\eps_i)}$ is uniformly bounded in the
$C^2$-topology, 
\item[(K2)]
$K_i |_{D(r_H)} \to H|_{D(r_H)}$ in the $C^2$-topology.
\end{itemize}
Note that $\pi \colon T^*M \to M$ is a Riemannian submersion with
respect to the Riemannian metrics $g^*$ and $g$. Applying
\eqref{e:nC} to $K_i$ we therefore find
\begin{equation}  \label{e1}
\mu_{g^*} \left( \gf_{K_i}^m \left( D_{q_0} (r_H+\eps_i )
\right)\right)               
\,\ge\, 
\left( c_1 \gb r_H \right) m \mu_g \left( W (q_0,K_i^m) \cap W (q_0,G_\gb^m)
\right) .
\end{equation}
Since $\eps_i \to 0$, we have
\[
\lim_{i \to \infty} \mu_{g^*} \bigl( \gf_{K_i}^m \left(
A(\eps_i) \right) \bigr) 
\,=\, 
\lim_{i \to \infty} \mu_{g^*} \bigl( \gf_{G_\gb}^m \left( A(\eps_i) \right) \bigr)
\,=\, 0 ,
\]
so that, together with Lemma~\ref{l:full}, 
\begin{equation}  \label{e2}
\lim_{i \to \infty} \mu_g \bigl( W (q_0,K_i^m) \cap W
(q_0,G_\gb^m) \bigr) 
\,=\,
\lim_{i \to \infty} \mu_g \bigl( V (q_0,K_i^m) \cap V (q_0,G_\gb^m) \bigr) 
\,=\, \mu_g(M) .
\end{equation}
Moreover, $\eps_i \to 0$ and (K1) imply
\begin{equation}  \label{e3}
\lim_{i \to \infty} \mu_{g^*} \Bigl( \gf_{K_i}^m 
\bigl( D_{q_0} (r_H+\eps_i) \setminus D_{q_0} (r_H) \bigr) \Bigr)
\,=\, 0 ,
\end{equation}
and (K2) implies
\begin{equation}  \label{e4}
\lim_{i \to \infty} \mu_{g^*} \bigl( \gf_{K_i}^m \left( D_{q_0}
(r_H) \right) \bigr) 
\,=\,
\mu_{g^*} \bigl( \gf_H^m \left( D_{q_0} (r_H) \right) \bigr) .
\end{equation}
Using \eqref{e4}, \eqref{e3}, \eqref{e1} and \eqref{e2} we find that
\begin{equation}  \label{e:last}
\mu_{g^*} \bigl( \gf_H^m \left( D_{q_0} (r_H) \right) \bigr) \,\ge\, 
\bigl( c_1 \gb r_H \mu_g(M) \bigr) m
\end{equation}
for $m$ large enough.
Since $q_0 \in M$ was arbitrary, Theorem~1\,(i) follows.

Assume now that $M$ is rationally elliptic but not simply
connected.
Then the universal cover $\widetilde{M}$ is rationally elliptic
and simply connected.
Let $pr \colon T^* \widetilde{M} \to T^* M$ be the projection
induced by the projection $\widetilde{M} \to M$,
and let $\tilde{g}^* = pr^* g^*$ and $\widetilde{H} = H \circ pr$.
If $H$ meets assumption~\eqref{ge:r}, so does $H$.
For $q_0 \in M$ we choose $\tilde{q}_0 \in \widetilde{M}$ over $q_0$ 
and notice that for each $m$
the projection $pr$ maps $\gf_{\widetilde{H}}^m \left(
D_{\tilde{q}_0} (r_H) \right)$ isometrically to $\gf_H^m \left(
D_{q_0} (r_H) \right)$.
Together with \eqref{e:last} we conclude that
\[
\mu_{g^*} \bigl( \gf_H^m \left( D_{q_0} (r_H) \right) \bigr)
\,=\,
\mu_{\tilde{g}^*} \left( \gf_{\widetilde{H}}^m 
\left( D_{\tilde{q}_0} (r_H) \right) \right) 
\,\ge\, 
\left( c_1 \gb r_H \mu_{\tilde{g}} \bigl( \widetilde{M} \bigr) \right) m ,
\]
so that Theorem~1 also follows for $(M,g)$.

\b
\ni
{\bf Step~2. From Floer homology to the homology of the path space} 

\m
\ni
Viterbo \cite{V1, V2} was the first to notice that the Floer homology for
periodic orbits of $T^*M$ is isomorphic to
the singular homology of the loop space of $M$. 
Different proofs were found by Salamon-Weber~\cite{SW} 
and Abbondandolo-Schwarz~\cite{AS}.
The work \cite{AS} also establishes a relative version of this
result: The Floer homology for Lagrangian intersections of
$T^*M$ is isomorphic to the singular homology of the based loop
space of $M$.
It is this version that we shall take advantage of.

We abbreviate $\rho = r_H + 2\eps$.
In order to estimate the dimension of $\HF_* \left( q_0,q_1,
G_\gb^m,\rho \right)$ from below, we first describe these
$\QQ$-vector spaces in a somewhat different way.
Set
\[
\cp^{m\gb \rho^2} \left( q_0,q_1,G_\gb^m \right)
\,:=\, \left\{ x \in \cp \left( q_0,q_1,G_\gb^m \right)
\mid \ca_{G_\gb^m} (x) \le m\gb \rho^2 \right\} .
\] 
\begin{lemma}  \label{l:rho}
$\cp^{m\gb \rho^2} (q_0,q_1,G_\gb^m) = \cp (q_0,q_1,G_\gb^m,\rho)$.
\end{lemma}

\proof
Assume that $x \in \cp \left( q_0,q_1,G_\gb^m \right)$.
For $t_0 \in [0,1]$ we choose geodesic normal coordinates $q$ 
near $\pi(x(t_0))$ in $M$. With respect to these coordinates the
equation $\dot{x}=X_{G_\gb^m}(x)$ at $t_0$ reads
\begin{eqnarray*}
\dot{p}(t_0) &=& 0,\\
\dot{q}(t_0) &=& 2m\gb p(t_0).
\end{eqnarray*}
Therefore,
$\gl \left( \dot{x}(t_0)\right) - G_\gb^m(x(t_0)) = 2 m\gb
|p|^2-m\gb|p|^2 = m\gb|p|^2=G_\gb^m \left(x(t_0)\right)$.
Since $G_\gb^m$ is autonomous, $G_\gb^m(x(t_0))$ does not depend
on $t_0$.
Integrating over $[0,1]$ we thus obtain
\[
\ca_{G_\gb^m} (x) \,=\, G_\gb^m (x).
\]
The lemma follows.
\proofend

We define the $k^{th}$ Floer chain group
$\CF^{m\gb \rho^2}_k \left( q_0,q_1,G_\gb^m \right)$ 
as the finite-dimensional $\QQ$-vector space freely
generated by the elements of 
$\cp^{m\gb \rho^2} \left( q_0,q_1,G_\gb^m \right)$ of index $k$.
Lemma~\ref{l:rho} yields
\begin{lemma}  \label{l:CF}
$\CF^{m\gb \rho^2}_* \left( q_0,q_1,G_\gb^m \right) = 
 \CF_* \left( q_0,q_1,G_\gb^m,\rho \right)$.
\end{lemma}

Denote by $\widehat{\cj}$ the set of families $\widehat{\bJ} =
\bigl\{ \widehat{J}_t \bigr\}$ of almost complex structures on
$T^*M$
such that $\widehat{\bJ} |_{D(\rho)} \in \cj$ 
and $\widehat{\bJ}$ is invariant under the flow of $Y$ on 
$T^*M \setminus D(r_H+\eps)$.
For $\widehat \bJ \in \widehat\cj$, 
for smooth maps $u \colon S \to T^*M$, 
and for $x^\pm \in \cp^{m\gb \rho^2} \left( q_0,q_1,G_\gb^m \right)$ 
consider Floer's equation
\begin{eqnarray}  \label{floer2}
 \left\{       
  \begin{array}{lcr}
    \pp_s u+\widehat{J}_t(u) \left( \pp_t u - X_{G_\gb^m}(u)
                                \right) = 0,  \\ [0.2em] 
    u(s,j) \in T^*_{q_j}M,\,\, j=0,1 , \\ [0.2em]
    \displaystyle{\lim_{s \to \pm \infty}} u(s,t) \,=\, x^\pm (t) 
             \text{ uniformly in } t .
  \end{array}
 \right.
\end{eqnarray}
We denote the set of solutions of~\eqref{floer2} by 
$\cm^{m\gb \rho^2} \left( x^-,x^+,G_\gb^m \right)$.
Lemmata~\ref{l:rho} and \ref{l:convex} imply
\begin{lemma}  \label{l:cm}
$\cm^{m\gb \rho^2} \left( x^-,x^+,G_\gb^m \right) = 
 \cm \left( x^-,x^+,G_\gb^m,\rho \right)$.
\end{lemma}
A standard argument shows that $\ca_{G_\gb^m} (x^-) \ge
\ca_{G_\gb^m} (x^+)$ for each 
$u \in \cm^{m\gb \rho^2} \left( x^-,x^+,G_\gb^m \right)$.
For generic $\widehat{\bJ} \in \widehat{\cj}$ 
the usual definition of the Floer boundary operator therefore
yields boundary operators
\[
\pp_k \bigl( \widehat{\bJ} \bigr) \colon
\CF_k^{m\gb \rho^2} \left( q_0,q_1,G_\gb^m \right) \,\to\,
\CF_{k-1}^{m\gb \rho^2} \left( q_0,q_1,G_\gb^m \right).
\]
Their homology groups $\HF_k^{m\gb\rho^2} \left( q_0,q_1,G_\gb^m
\right)$ do not depend on $\widehat \bJ$.
Lemmata~\ref{l:CF} and \ref{l:cm} imply that 
$\pp_k \bigl( \widehat{\bJ} \bigr) = \pp_k \bigl( \bJ \bigr)$, 
whence
\begin{proposition}  \label{p:floer=}
$\HF_* \left( q_0,q_1,G_\gb^m,\rho \right)
\,\cong \, 
\HF_*^{m\gb\rho^2} \left( q_0,q_1,G_\gb^m\right)$.
\end{proposition}

For $q_0, q_1 \in M$ let $\Omega^1 \left( M,q_0,q_1\right)$
be the space of all paths $q \colon [0,1] \to M$ of Sobolev
class $W^{1,2}$ such that $q(0)=q_0$ and $q(1)=q_1$.
Again, this space has a canonical Hilbert manifold structure.
The energy functional $\ce \colon \Omega^1 \left( M,q_0,q_1\right) \to \RR$
is defined as
\[
\ce (q) \,=\, \frac{1}{2} \int_0^1 \left| \dot{q}(t) \right|^2 dt.
\]
For $a>0$ we consider the sublevel sets 
\[
\ce^a(q_0,q_1) \,:=\, \left\{ q \in \Omega^1 \left( M,q_0,q_1\right) \mid
\ce(q) \le a \right\} .
\] 
In the following, $\H_*$ denotes singular homology with rational coefficients.
\begin{proposition} \label{p:AS} 
$\HF_*^{m\gb\rho^2} \left( q_0,q_1,G_\gb^m\right)
\,\cong \, 
\H_* \left( \ce^{2\left( \gb \rho m \right)^2} (q_0,q_1) \right)$.
\end{proposition}

\proof
Let $L \colon TM \to \RR$ be the Legendre transform of
$G_\gb^m$.
Applying Theorem~3.1 of \cite{AS} to $G_\gb^m$ and $L$, we
obtain
\[
\HF_*^{m\gb\rho^2} \left( q_0,q_1,G_\gb^m\right)
\,\cong \, 
\H_* \left( \left\{ q \in \Omega^1 (M,q_0,q_1) \mid \int_0^1 L
\left( q(t), \dot{q} (t) \right) dt \le m \gb \rho^2 \right\} \right) .
\]
Notice now that $L(q,v) = \frac{1}{4m\gb} \left|v\right|^2$.
The set $\{ \dots \}$ on the right hand side therefore equals
$\ce^{2\left( \gb \rho m \right)^2} (q_0,q_1)$, and so Proposition~\ref{p:AS}
follows.
\proofend 
        
\begin{remark*}
{\rm
In \cite{AS}, Abbondandolo and Schwarz work with almost complex
structures which are close to the almost complex structure 
interchanging the horizontal and vertical tangent bundles
of $\bigl( T^*M,g^* \bigr)$.
They need to work with such almost complex structures in order
to prove a subtle $L^\infty$-estimate for solutions of Floer's
equation which is crucial for obtaining their general result.
For the special functions $G_\gb^m$ appearing in our situation,
one can work with convex almost complex structures $\widehat \bJ
\in \widehat \cj$ and does not need their $L^\infty$-estimate.
\diam
}
\end{remark*}

Propositions~\ref{p:floer=} and \ref{p:AS} yield
\begin{proposition}  \label{p:floerlength}
$\HF_* \left( q_0,q_1,G_\gb^m,\rho \right)
\,\cong \, 
\H_* \bigl( \ce^{2\left( \gb \rho m \right)^2} (q_0,q_1) \bigr)$.
\end{proposition}

\ni
{\bf Step~3. Lower estimates for 
$\dim \H_* \bigl( \ce^{2 \left( \gb \rho m \right)^2} (q_0,q_1) \bigr)$.}  

\m
\ni
The length functional $\cl \colon \Omega^1 \left( M,q_0,q_1\right) \to \RR$
is defined as
\[
\cl (q) \,=\, \int_0^1 \left| \dot{q}(t) \right| dt.
\]
For $a>0$ we consider the sublevel sets 
\[
\cl^a(q_0,q_1) \,:=\, \left\{ q \in \Omega^1 \left( M,q_0,q_1\right) \mid
\cl(q) \le a \right\} .
\] 

\s
\ni
{\it Proof of Proposition~\ref{p:floer}\:(i).}
Throughout the proof of Proposition~\ref{p:floer}\:(i) we assume
that $(M,g)$ be a simply connected closed Riemannian manifold. 

\begin{lemma}  \label{l:gromov}
There exists a constant $C_G > 0$ depending only on $(M,g)$
such that each element of $\H_j \bigl( \Omega^1 (M,q_0,q_1)
\bigr)$ can be represented by a cycle in 
$\ce^{\frac 12 \left( C_G j \right)^2} (q_0,q_1)$.
In particular,
\[
\dim \H_* \bigl( \ce^{\frac 12 \left( C_G m \right) ^2}
(q_0,q_1) \bigr) 
\,\ge\,
\sum_{j=0}^m \dim \H_j \left( \Omega^1 (M,q_0,q_1) \right)
\quad\; \text{ for all } m .
\]
\end{lemma}

\proof
According to a result of Gromov, there exists a constant $C_G>0$ depending
only on $(M,g)$
such that each element of $\H_j \bigl( \Omega^1 (M,q_0,q_1)
\bigr)$ can be represented by  
a cycle lying in $\cl^{(C_G-2) j} (q_0,q_1)$.
Gromov's original proof of this result in \cite{Gr78} is very short.
Detailed proofs can be found in \cite{Pa1} and \cite[Chapter 7A]{Gr}.
Let $\Delta^j$ be the $j$-dimensional standard simplex,
and let $h \colon \Delta^j \to \cl^{(C_G-2) j} (q_0,q_1)$ be an
integral cycle.
By suitably reparametrizing each path $h(s)$ near $t=0$ and
$t=1$ and then smoothing each path with the same heat kernel we
obtain a homotopic and hence homologous cycle
$h_1 \colon \Delta^j \to \cl^{(C_G-1) j} (q_0,q_1)$
consisting of smooth paths.
We identify $h_1$ with the map
\[
\Delta^j \times [0,1] \to M, \quad\,  
(s,t) \mapsto h_1(s,t) := \left( h_1(s) \right) (t) .
\]
Endow the manifold $M \times [0,1]$ with the product Riemannian
metric, and set $\tilde{q}_0 = (q_0,0)$ and $\tilde{q}_1 =
(q_1,1)$.
We lift $h_1$ to the cycle 
$\tilde{h}_1 \colon \Delta^j \to \Omega^1 \left( M \times [0,1],
\tilde{q}_0, \tilde{q}_1 \right)$ defined by
$\tilde{h}_1(s,t) = \left( h_1(s,t),t \right)$.
This cycle consists of smooth paths whose tangent vectors do not vanish.
For each $s$ let $\tilde{h}_1 \left( \gs (s) \right)$
be the reparametrization of $\tilde{h}_1(s)$ proportional to arc length.
The homotopy 
$H \colon [0,1] \times \Delta^j \to 
\Omega^1 \left( M \times [0,1],\tilde{q}_0,\tilde{q}_1 \right)$ 
defined by 
\[
\bigl( H (\tau,s) \bigr) (t) \,=\, 
\tilde{h}_1 \bigl( s, (1-\tau) t + \tau \gs (s) \bigr)
\]
shows that $\tilde{h}_1$ is homologous to the cycle $\tilde{h}_2
(s) := H (1,s)$.
Its projection $h_2$ to $\Omega^1 (M,q_0,q_1)$ is homologous to
$h_1$ and lies in $\cl^{(C_G-1)j}(q_0,q_1)$.
Since for each $s$ the path $\tilde{h}_2(s)$ is parametrized
proportional to arc length, we conclude that
\begin{eqnarray*}
\ce \bigl( h_2(s) \bigr) 
  \,\le\, \ce \bigl( \tilde{h}_2(s) \bigr)
  &=& \tfrac 12 \left( \cl \bigl( \tilde{h}_2 (s) \bigr) \right)^2 \\
  &=& \tfrac 12 \left\{ \left( \cl \bigl( h_2 (s)
                                   \bigr) \right)^2 +1 \right\} 
\,\le\, \tfrac 12 \left( C_G-1 \right)^2 j^2 + \tfrac 12
\,\le\, \tfrac 12 \left( C_G j \right)^2 
\end{eqnarray*}
for each $s$, so that indeed $h_2 \subset \ce^{\frac 12 (C_Gj)^2}(q_0,q_1)$.
%
%
\proofend

Let $\Omega (M, q_0,q_1)$ be the space of {\it continuous}\, path
$q \colon [0,1] \to M$ from $q_0$ to $q_1$
endowed with the compact open topology.
According to \cite[Chapter~17]{Mi1} or \cite[Theorem~1.2.10]{Klin-lectures},
the inclusion $\Omega^1 (M, q_0,q_1) \to \Omega (M, q_0,q_1)$ 
is a homotopy equivalence.
The homotopy type of these spaces does not depend on $q_0,q_1$
and is denoted $\Omega (M)$.
Together with Proposition~\ref{p:floerlength} and
Lemma~\ref{l:gromov} we find
\begin{equation}  \label{e:CG}
\dim \HF_* \left( q_0,q_1,G_\gb^m,\rho \right) 
\,\ge\, 
\sum_{j=0}^{\big\lfloor \frac{2}{C_G} \gb \rho m \big\rfloor} 
\dim \bigl( \H_j (\Omega(M) ) \bigr) 
\quad\; \text{ for all } m .
\end{equation} 

Since $M$ is a simply connected and closed manifold of dimension
$d$, 
Proposition~11 on page~483 of
Serre's seminal work~\cite{Se} guarantees that 
for every integer $i \ge 0$ there exists an integer $j \in
\left\{ 1, \dots, d-1 \right\}$ 
such that
$\H_{i+j} \bigl( \Omega (M) \bigr) \neq 0$. 
Since $\H_0 \bigl( \Omega (M) \bigr) \neq 0$, we find that
\[
\sum_{j=0}^m \dim \H_j \bigl( \Omega (M) \bigr)
\,\ge\,
1+ \bigg\lfloor \frac m d \bigg\rfloor \,\ge\, \frac m d
\quad \; \text{ for all } m .   
\]
Setting $c_1 := 1/(d C_G)$ we conclude together with the
estimate~\eqref{e:CG} that 
\[
\dim \HF_* \left( q_0,q_1,G_\gb^m,\rho \right) 
\,\ge\, 
c_1 \gb \rho m 
\quad \; \text{ for all large enough } m ,
\]
proving the first claim of Proposition~\ref{p:floer}\:(i).

Assume now that $(M,g)$ is rationally hyperbolic.
Then there exists $C>1$ such that
\[
\sum_{j=0}^m \dim \pi_j (M) \otimes \QQ \,\ge\, C^m
\quad \; \text{ for all large enough } m .
\]
According to \cite{MM}, 
$\dim \pi_{j+1} (M) \otimes \QQ \le 
\dim \H_j \left( \Omega (M) \right)$ 
for all $j \ge 0$, so that
\begin{equation}  \label{e:MM}
\sum_{j=0}^m \dim \H_j \bigl( \Omega (M) \bigr) \,\ge\, C^m
\quad \; \text{ for all large enough } m .
\end{equation}
With $C_1 := \frac {1}{C_G} \log C$,
the estimates~\eqref{e:CG} and \eqref{e:MM} yield the desired estimate
\[
\dim \HF_* \left( q_0,q_1,G_\gb^m,\rho \right) 
\,\ge\, 
e^{C_1 \gb \rho m} \quad \; \text{ for all large enough } m .
\]

\begin{remark}  \label{r:rational}
{\rm
In order to have the estimate \eqref{e:MM} at hand, it is
important that we can work with rational coefficients.
Indeed, let $M$ be a simply connected closed manifold such that 
the sequence
$\sum_{j=0}^m \dim \H_j \left( \Omega (M); \ZZ_2 \right)$
grows faster than every polynomial in $m$.
Then it is only known that this sequence grows faster than
$C^{\sqrt{m}}$ for some constant $C>1$, see~\cite{FHT93}.
}
\end{remark}

\m
\ni
{\it Proof of Proposition~\ref{p:floer}\:(ii).}
Assume that $\pi_1(M)$ is infinite.
We first notice that
\begin{equation}  \label{e:first}
\dim \H_* \bigl( \ce^{2 \left( \gb \rho m \right)^2} (q_0,q_1) \bigr) 
\,\ge\,
\dim \H_0 \bigl( \ce^{2 \left(\gb \rho m \right)^2} (q_0,q_1) \bigr) \\
\,=\,
\# \pi_0 \bigl( \ce^{2 \left( \gb \rho m \right)^2} (q_0,q_1) \bigr) .
\end{equation}
For $a>0$ denote by $\Pi^a_\ce (q_0,q_1)$ 
(resp.~$\Pi^a_\cl (q_0,q_1)$)  
the set of those homotopy classes of $W^{1,2}$-path $q \colon [0,1] \to M$ 
from $q_0$ to $q_1$ which can be represented by a path
of energy (resp.~length) at most $a$.
Then 
\begin{equation}  \label{e:second}
\# \pi_0 \bigl( \ce^{2 \left( \gb \rho m \right)^2} (q_0,q_1) \bigr) 
\,\ge\,
\# \Pi_\ce^{2 \left( \gb \rho m \right)^2} (q_0,q_1)
\,=\, 
\# \Pi_\cl^{2 \gb \rho m} (q_0,q_1) .
\end{equation}
Choose a smooth path $h$ from $q_0$ to $q_1$ with $\length (h) \le \delta
:= \dia (M,g)$,
and assume that $m$ is so large that $2 \gb \rho m > \delta$.
Since the map
\[
\Pi_\cl^{2\gb\rho m-\gd} (q_0,q_0) \to \Pi_\cl^{2\gb \rho m}(q_0,q_1) ,
\quad \;
[\go] \mapsto [h \circ \go] ,
\]
is injective, we have that
\begin{equation}  \label{e:third}
\# \Pi_\cl^{2\gb \rho m} (q_0,q_1)
\,\ge\, 
\# \Pi_\cl^{2\gb \rho m - \gd} (q_0,q_0).
\end{equation}
Let now $S = \{ h_1, \dots, h_{\#S} \}$ be a generating set of
$\pi_1(M)$.
In view of the definition of $l(S,g)$ in Section~\ref{s:gr}
we can represent each $h_j$ by a smooth loop based at $q_0$ 
of length no more than $l(S,g)$.
In view of the triangle inequality and the definition of the growth
function $\gamma_S$ we finally obtain
\begin{eqnarray}  \label{e:fourth}
\# \Pi_\cl^{2\gb \rho m - \gd} (q_0,q_0)
&\ge&
\gg_S \left( \lfloor (2\gb \rho m - \gd)/ l(S,g) \rfloor \right)\\
&\ge&
\gg_S \left( \lfloor (2\gb r_H) m / l(S,g) \rfloor \right) \notag
\end{eqnarray}
for $m$ large enough.
Proposition~\ref{p:floerlength} and the estimates 
\eqref{e:first}, \eqref{e:second},  
\eqref{e:third}, \eqref{e:fourth} yield
\[
\dim \HF_* (q_0,q_1, G_\gb^m,\rho) \,\ge\, 
\gg_S \bigl( \lfloor (2 \gb r_H) m / l(S,g) \rfloor \bigr) , 
\]
and so Proposition~\ref{p:floer}\:(ii) follows.
\proofend

\subsection{Proof of Corollary~1}

For $M=S^1$ the claim follows from an elementary topological
argument, see \cite{FS:GAFA}.
We can therefore assume that $(M,g)$ is a $P$-manifold of
dimension $d \ge 2$.
This has two consequences:
First, let $\Ham^c \left( T^*M \right)$ be the group of
diffeomorphisms of $T^*M$ generated by compactly supported
Hamiltonians $H \colon [0,1] \times T^*M \ra \RR$,
and let $\Symp_0^c \left( T^*M \right)$ be the group of 
diffeomorphisms of $T^*M$ which are isotopic to the identity
through a family of symplectomorphisms supported in a 
compact subset of $T^*M$.
Then
\[
\Ham^c \left( T^*M \right) \,=\, \Symp_0^c \left( T^*M \right) 
\]
(see \cite[Lemma~2.18]{FS:GAFA}). 
Moreover, $M$ is rationally elliptic (see~Proposition~\ref{p:Pr} below).
Let now $\psi \in \Symp^c \left( T^*M \right)$ be such that
$[ \psi ] = \left[ \gt^m \right] \in \pi_0 \left( \Symp^c \left( T^*M
\right) \right)$ for a twist $\gt = \gt_f$ on $T^*M$ and some $m \in \ZZ
\setminus \{ 0 \}$.
Then $\psi \gt^{-m} \in \Symp^c_0 \left( T^*M \right) = \Ham^c
\left( T^*M \right)$, so that we find a compactly supported
Hamiltonian function $H \colon [0,1] \times T^*M \to \RR$ with
$\gf_H = \gt^{-m} \psi$.
Then $\psi = \gt^m \gf_H = \gf_{mf} \gf_H = \gf_K$, where
\[
K(t,q,p) \,=\, m f \left( |p| \right) + H \left( t, \gt^{-m}_f (q,p) \right) .
\]
Choose $r_k \ge 1$ so large that $H$ is supported in $[0,1]
\times \left\{ |p| \le r_K \right\}$.
Since $\gt_f^{-m}$ preserves the levels $\left\{ |p| = \const
\right\}$, we then have 
$K(t,q,p) \,=\, m f \left( |p| \right)$ for $|p| \ge r_k$.
Since $f'(r_K) =1$ and $m \neq 0$, Theorem~1\:(i) applies and
yields $\check{s}_{\fibre} \left( \gf_K;r_K\right) \ge 1$, as claimed.
\proofend

\section{More on $P$-manifolds}  \label{s:P}

\ni
Much information on $P$-manifolds can be found in the book
\cite{B} and in Section~10.10 of \cite{Be}.
In this section we give a few additional results.

\begin{proposition}  \label{p:Pr}
Let $(M,g)$ be a $P$-manifold of dimension at least~$2$. 
Then $M$ is rationally elliptic. 
\end{proposition}

\proof  
By the Bott--Samelson theorem of B\'erard Bergery \cite[Theorem~7.37]{B}
the fundamental group $\pi_1(M)$ is finite.

\s
\ni
{\it First argument:}
Since $\pi_1(M)$ is finite, the universal covering 
$\bigl( \widetilde{M},\tilde{g} \bigr)$ is also a $P$-manifold,
and the rational homotopy groups of $M$ and $\widetilde{M}$ are the same.
We can thus assume that $M$ is simply connected.
Recall that $\dim \pi_{j+1} (M) \otimes \QQ \le \dim \H_j \left(
\Omega (M) \right)$ for all $j \ge 0$.
It therefore suffices to show that the numbers 
$\dim \H_j \left( \Omega M \right)$ are uniformly bounded.
Recall that we scaled $g$ such that all unit-speed geodesics
have minimal period $1$.
For such a geodesic $\gg \colon \RR \to M$ and $t>0$
we let $\ind \gg (t)$ be the number of linearly independent Jacobi
fields along $\gg (s)$, $s \in [0,t]$, 
which vanish at $\gg (0)$ and $\gg (t)$.
If $\ind \gg (t) >0$, then $\gg (t)$ is said to be conjugate to 
$\gg (0)$ along $\gg$.
The index of $\gg |_{[0,a]}$ defined as
\[
\ind \gg |_{[0,a]} \,=\, \sum_{t \in ]0,a[} \ind \gg (t) 
\]
is a finite number, 
and according to \cite[1.98 and 7.25]{B}
the number $k = \ind \gg |_{[0,1]}$ is the same for all 
unit-speed geodesics $\gg$ on $(M,g)$.
Fix now $q_0 \in M$ and choose $q_1$ which is not
conjugate to $q_0$. Then there are only finitely many geodesic segments
$\gg_1, \dots, \gg_n$ from $q_0$ to $q_1$ of length smaller
than~$1$, see \cite[7.41]{B}.
Let $k_1, \dots, k_n$ be their indices.
The energy functional on $\Omega^1 (M;q_0,q_1)$ is Morse with indices 
$k_i+l(d-1+k)$, where $1 \le i \le n$ and $l = 0,1,2, \dots$.
Since $d-1+k \ge 1$, we conclude that
$\dim \H_j \left( \Omega (M) \right) = \dim \H_j \left( \Omega
(M;q_0,q_1) \right) \le n$. 

\s
\ni
{\it Second argument:}
Recall from Discussion~2\:(i) that an elementary
computation shows that $\check{s}_{\fibre} ( \gt_f ) =1$ for every
twist on $T^*M$. The proposition thus follows from Theorem~1\:(i).
\proofend

The main statement in the Bott--Samelson theorem for
$P$-manifolds \cite[Theorem~7.37]{B} is that the rational
cohomology ring of such a manifold has only one generator.
Comparing with the lists given in \ref{ss:rational} we find 

\begin{proposition}
Assume that $(M,g)$ is a $P$-manifold of dimension~$d$, and
denote by $\widetilde{M}$ its universal covering.

\s
If $d=3$, then $\widetilde{M}$ is diffeomorphic to $S^3$.

\s
If $d=4$, then $\widetilde{M}$ is homeomorphic to $S^4$ or $\CP^2$.

\s
If $d=5$, then $\widetilde{M}$ is diffeomorphic to $S^5$ or $\SU
(3) / \SO (3)$.
\end{proposition}

\begin{remark}
{\rm
We do not know whether $\SU (3) / \SO (3)$ carries a \text{$P$-metric}.
A $P$-metric is said to be an $SC$-metric if all closed geodesics 
are embedded circles of equal length.
It follows from the Bott--Samelson theorem for $SC$-manifolds
\cite[Theorem~7.23]{B}
and from $H^2\left( \SU (3) / \SO (3) ;\ZZ \right) =
\ZZ_2$
that $\SU (3) / \SO (3)$ cannot carry an $SC$-metric.
\diam
}
\end{remark}

\section{Outlook}

\ni
The conceptual point of view of this paper was to look at entropy-type
quantities which are well understood for geodesic flows,
and to establish the lower bounds for these quantities known for
geodesic flows for arbitrary 
classical Hamiltonian systems by interpreting these
quantities in Floer homological terms and by using the deformation
invariance of Floer homology.
We were able to do this for Hamiltonians meeting \eqref{ge:r} by
using the Convexity Lemma~\ref{l:convex}.
Already the autonomous Hamiltonians 
$H(q,p) = \frac 12 \left| p - A(q) \right|^2 +V(q)$ modelling
the dynamics of a particle in a magnetic potential $A(q)$ and a scalar
potential $V(q)$ do not meet \eqref{ge:r},
and the Convexity Lemma fails for such Hamiltonians.
In \cite{MaS} the $L^\infty$-estimate for solutions of Floer's
equation from \cite{AS} is used to extend the results of this
paper to Hamiltonians which are autonomous and convex above some
energy level.

%


\enddocument